# Distributed Load Scheduling in Residential Neighborhoods for Coordinated Operation of Multiple Home Energy Management Systems

Babak Jeddi*, Yateendra Mishra, Gerard Ledwich

*School of Electrical Engineering and Robotics, Queensland University of Technology, Brisbane, Australia*

**Abstract-** Recently, home energy management systems (HEMS) are gaining more popularity enabling customers to minimize their electricity bill under time-varying electricity prices. Although they offer a promising solution for better energy management in smart grids, the uncoordinated and autonomous operation of HEMS may lead to some operational problems at the grid level. This paper aims to develop a coordinated framework for the operation of multiple HEMS in a residential neighborhood based on the optimal and secure operation of the grid. In the proposed framework customers cooperate to optimize energy consumption at the neighborhood level and prevent any grid operational constraints violation. A new price-based *global* and *individualized* incentives are proposed for customers to respond and adjust loads. The individual customers are rewarded for their cooperation and the network operator benefits by eliminating rebounding network peaks. The alternating direction method of multipliers (ADMM) technique is used to implement coordinated load scheduling in a distributed manner reducing the computational burden and ensure customer privacy. Simulation results demonstrate the efficacy of the proposed method in maintaining nominal network conditions while ensuring benefits for individual customers as well as grid operators.

## 1    Introduction

Nowadays, rooftop solar panels are very common and affordable worldwide, specifically in Australia where 1 in every 5 households has rooftop PV [1]. Also, solar customers can make more profit by installing battery energy storage systems (BESS) which enable them to store energy produced by PV instead of wasting it into the grid at a small rate [2]. On the other side of the network, utility companies and network operators try to maintain a balance between the demand and generation within the network by launching demand response (DR) programs to engage customers in modifying their electricity usage in response to the network operating conditions [3, 4]. The implementation of such mechanisms requires two-way communication between the utility companies and the customers facilitated by the arrival of smart metering devices and communication infrastructures [5].

Existing DR methods are separated into either direct load control (DLC) techniques or price-based techniques [6]. In the DLC mechanism, the network operator can directly control the operation of certain devices of customers (e.g., air-conditioning systems). In return, customers would receive financial incentives. In the price-based DR mechanism, whereas, the network operator offers time-variant price signals to the customers anticipating that they respond to the prices and adjust their consumptions accordingly [7]. There are various pricing methods such as real-time pricing (RTP), time of use (ToU) pricing, and critical peak pricing (CPP) [6]. To respond to the price signals, customers need to have an intelligent energy management system to assist them with making reasonable consumption decisions and simultaneously get the most out of their solar system and BESS. In this regard, the popularity of home energy management systems (HEMS) is growing rapidly [8, 9]. A HEMS determines the best scheduling for the operation of household devices based on the price signals given the operational constraints of devices so that the electricity bill is minimized.

---

* Corresponding author at School of Electrical Eng. and Robotics, Queensland University of Technology, Brisbane City, QLD 4000, (email: babak.jeddi@hdr.qut.edu.au)

Recently, numerous studies have focused on developing new HEMS models to save more money for households using powerful optimization and management techniques [2]. Although these studies provide significant financial benefits for individual customers, most of them overlook the potential undesirable impacts on the distribution grid. When a group of HEMS works independently under the same time-dependent pricing scheme, every household will simultaneously schedule the load to the low-price periods, and consequently, a "rebound peak" may occur at low price hours [10, 11]. Furthermore, the local transformer serving that residential area could experience an "overload" which may lead to failures and accelerated aging [12]. One possible strategy to tackle these problems is to coordinate the operation of multiple HEMS such that the aggregated load profile of the feeder is improved [13, 14]. Such coordinated energy management can be implemented in either a centralized or distributed fashion. In centralized mode, one central controller is controlling the operation of all HEMS, while in the distributed mode, each HEMS is being controlled independently but communicating with the central controller or others in the neighborhood area [15, 16]. While the centralized model gives an effective result; it suffers from the heavy computational burden and requires all detailed information of the customers. The distributed model, whereas, requires a minimum amount of information exchange and can be implemented efficiently making it more scalable and hence is a  preferred option [17].

## 2   Literature Review and Contributions

A centralized coordinated energy management platform is proposed to mitigate the unfair usage of the distribution transformer capacity and minimize the total energy procurement cost of the neighborhood [18]. In [19], a centralized DR model is developed for a large neighborhood including residential, commercial and industrial customers aiming to control the operation of devices across the network for minimizing electricity payment and improving the aggregated load profile of the network. The model aims to create a desirable load curve while minimizing the energy cost and peak load demand using a heuristic-based Evolutionary Algorithm. In [20], an aggregator-based DR program is developed in which the aggregator maximizes its profit by selling the capacity aggregated from customer smart appliances. In real-time, each customer has the option to buy electricity from the utility at the real-time price or the aggregator at a price named customer incentive price. The aggregator collects settings from all end-users and calculates optimal setpoints using a genetic algorithm-based method for all houses. A coordinated energy management problem is formulated in [21] as a bi-level optimization model in which the upper sub-problem aims to adjust the networks load profile and the lower sub-problems minimize each customers' electricity expenditures. To make it simple, the bi-level problem is converted to a single-level optimization problem in a decentralized fashion; first, homes schedule their electricity consumption based on the price announced by the central controller and send the results back. Then, the central controller computes the aggregated load curve and sends it with a desired aggregated curve back to the customers. Then, customers schedule their electricity consumption again to make the aggregated curve as much as close to the desired one. This method was effective in benefiting distribution system operators in terms of load profile improvement; however, it suffers from convergence and robustness problems. In Addition, customers are not financially incentivized for their cooperation in networks load profile improvement.

The load management problem formulated in [22] aims at maximizing the total utility of all customers to reduce peak load demand in the network. A distributed dynamic programming approach is employed to solve the problem in a distributed way in which each customer is represented by an agent who can exchange information with its neighboring agents. During load management, an agent first receives the information of load settings and incentives by the system operator. Then, the agent participates in the optimization process in cooperation with its neighboring agents to obtain a solution. In [23], a collaborative energy management algorithm is developed to improve the aggregated consumption profile and decrease the cost of real-time power balancing in the neighborhood. Customers are connected to a central controller and they interact with each other by exchanging messages. The problem is formulated as a multi-stage stochastic optimization model and several approximations and decomposition methods applied to solve it in a decentralized fashion.

To sum up, some research work has been conducted in this area, however, there are still some aspects of the topic that has not been given enough investigations. Firstly, the majority of the existing coordinated energy management models have ignored including network constraints in the formulation. Network constraints in terms of line current limits are especially important as the number of HEMS is growing, and more unwanted consumption peaks are likely to happen. Secondly, appropriate financial incentives and keeping the comfort level of customers have not been fully considered in the literature previously. Customer's participation in any coordinated energy management scheme considerably depends on the financial motivation offered to them. Most importantly, they should be able to maintain their energy usage privacy and do not encounter any deterioration in energy service quality.

This paper, therefore, proposes a framework for the coordinated operation of multiple HEMS to minimize electricity payment for the individual customer and to keep the operational constraints of the network within desirable limits. The proposed coordinated load scheduling model is implemented in a distributed manner as it offers more flexibility and freedom to the customers and preserves their privacy. The main goal of the network operator is always to achieve a safe operational condition for the network and balance the total demand with the available supply. On the other side, individual customers always seek to minimize their electricity bills. These two objectives are considered in the proposed model and the electricity price is adjusted to incentivize customers to schedule their loads in a way that not only minimize their bill but also benefits the network. In other words, by participating in such a coordinated consumption program, the payment of customers does not increase compared to the case when they are operating independently, and they do not experience any degree of discomfort in operating their devices.

In the proposed framework, HEMS schedule the controllable devices over the day considering the inhabitant's desires in line with the network's operational conditions. Putting differently, each customer is not only concerned about its energy consumption but also cooperates with others to optimize energy consumption at the neighborhood level and prevents any limitation violation in the network's operational constraints, and accordingly, get compensations for responding to the network's issues.

Available techniques to implement the coordinated load scheduling model in a distributed manner include game theory [24], multi-agent systems [25], dual decomposition approach [26], bi-level optimization [13, 21], analytical target cascading [27] and column generation + Dantzig-Wolfe decomposition method [11]. Each of these methods has its advantages and drawbacks. For example, in the method of multi-agent systems, agents can be automatically adjusted to new environmental variations such as changes in the network structure with minimum modification needed in the algorithm structure. However, they require a large amount of communicating units making them costly for real-world implementation. For a detailed comparison for different distributed optimization techniques, please refer to [16, 17].

In this paper, the proposed distributed optimization formulation is based on the alternating direction method of multipliers (ADMM) method [28]. This method decomposes the large-scale problem into several small-scale problems and allows customers to solve their energy management problem independently in parallel with other customers with a minimum exchange of information to the network operator.

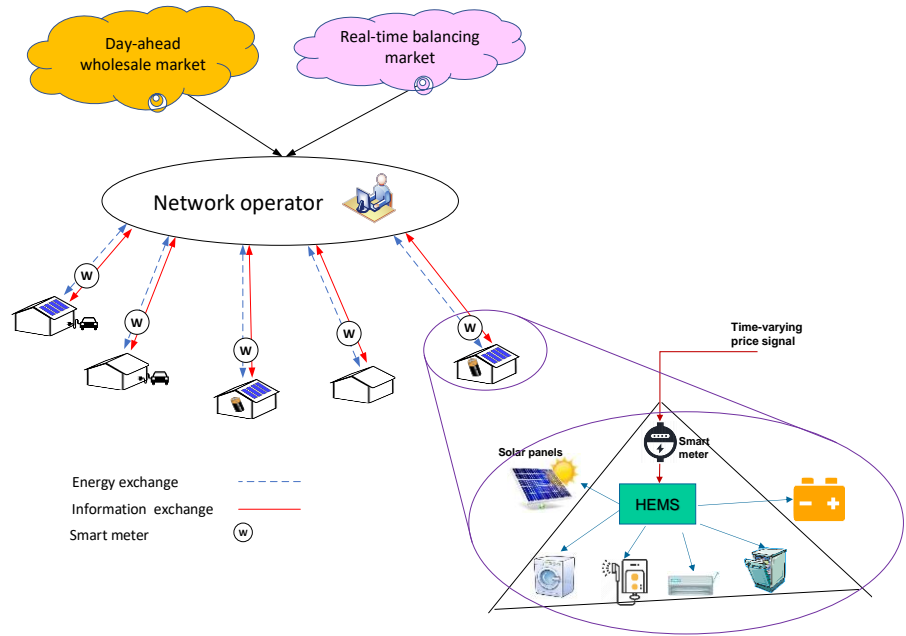

Figure 1: Network structure for implementing the proposed coordinated load scheduling model

In summary, the novelties of the proposed model can be summarized as follows:

1- A distributed load scheduling algorithm in a residential network considering network operational constraints is proposed.

2- A new mechanism for designing appropriate incentive for customers is formulated and *global* and *individualized* pricing structures are proposed.

3- Customers independently schedule the operation of their devices and only send their total consumption at each hour of the day to the central controller. They do not share any other details about their electricity decisions.

The rest of this paper is organized as follows: In Section 3, a general description of the proposed load scheduling framework is presented. Then, a mathematical model of customers and network operator is formulated in Section 4. The centralized and distributed optimization model along with the new pricing mechanisms are formulated in this section as well. Finally, simulation results are presented in Section 5 followed by conclusions and future work in Section 6.

## 3    Proposed Model Description and Assumptions

A test network connected to the upstream transmission grid through a step-down substation transformer supplying electricity to a certain number of residential customers ( Figure 1) is considered in this paper. An aggregator is responsible for supplying power to the distribution network. The aggregator can be a utility company who owns and operates the network, or it can be a separate third-party entity. If the aggregator is the distribution network operator (DNO), it oversees all the costs associated with the system. However, if the DNO is a separate third-party entity, it must reimburse the DNO for the damages on network components due to its interaction with customers.  In this paper, it is assumed that aggregator and DNO are not separate entities and there is only one entity that owns the network and is responsible for the operation of the network. Thus, the "aggregator", "network operator" and "utility company" terms are used interchangeably throughout the paper.

Each customer has installed a HEMS and there are a certain number of customers with solar PVs and BESS who can sell surplus power back into the network at a feed-in-tariff (FiT) price. The BESS can be charged from the PV generation or imported power from the network. Each home is equipped with an advanced metering infrastructure (AMI) to communicate directly with a central controller through a protected local area network (LAN). It is assumed that the installation cost of AMI and all other communication technologies have already been paid by the DNO.

A day-ahead pricing mechanism is utilized in which the price of electricity for the next day is decided on the previous day. The DNO buys energy from the upstream grid and sells it to the customers. It participates in the day-ahead wholesale market and a real-time balancing market. It needs to forecast the total consumption in the network for the next day and accordingly purchase power from the day-ahead market. Additionally, the DNO might own some distributed energy resources within the network and utilize them for suppling the network load. Then, it should determine dynamic base prices based on the wholesale market price and its financial requirements and release them to the customers. The DNO should also buy additional power from the real-time market for balancing supply and demand within the network. Furthermore, it has to pay the upstream grid for absorbing the extra power form the real-time market that is not being consumed [23]. Therefore, DNO aims that the aggregated consumption within the network matches the available supply, i.e., the power that already purchased from the day-ahead market. It is to be noted that electricity pricing and its associated variability is a complex process and is beyond the scope of this paper. The prices here are just indicators of the real-world electricity markets and tariffs.

Once customers receive the base prices, installed HEMS are activated automatically to schedule household load. As mentioned before, considerable power imbalances might be created as the result of the independent operation of HEMS leading to an increase in the cost of the DNO in the real-time market. Thus, the DNO wishes to initiate a DR program to adjust price signals incentivizing customers to change their consumption/generation to fit the needs of the network. The amount of reduction or rise in price is proportional to the gap that exists between the available supply and the actual load of the customers. Thus, customers try to cooperate in reshaping the network's load profile and coordinate their consumption aiding the DNO to manage the demand/supply balance and associated congestion in the network.

It is assumed that the DNO and all customers have already reached an agreement and all customers have agreed to participate in this program. Additionally, customers do not cut or reduce their daily consumption or generation of energy, rather just shift it in time during the day from the DNO. This is performed by changing the operation time of controllable appliances or changing the charging/discharging profile of the BESS. In this paper, all the uncertainties associated with household electricity consumption, price, and PV generation are ignored and a deterministic model is developed.

## 4    Mathematical Formulation

### 4.1    Day-ahead price setting

The DNO determines time-varying base prices for the next day based on the forecasted demand of the network considering its financial and technical preferences. Dynamic base prices are usually proportional to the available supply, i.e., they are high during the high-supply hours and low at low-supply periods. It is assumed the DNO calculates the day-ahead base price $\pi_t$ for each time slot $t$ of the day using (1) below [29, 30].

$$\pi_t = \varphi \frac{\partial F_t^{DHA}(g_t)}{\partial g_t} \tag{1}$$

where $g_t$ is the total forecasted demand of the network at time $t$; $F_t^{DHA}(.)$ is the cost function in the day-ahead market; and $\varphi > 1$ is the profit coefficient.

There are various forecasting techniques proposed in the literature to compute $g_t$ which are beyond the scope of this paper. Here, for simplicity, (2) is used in this paper.

$$g_t = \sum_{h=1}^{H} \widetilde{D}_{h,t} + \widetilde{Loss}_t - P_t^{DG} \tag{2}$$

where $H$ is the total number of homes connected to the network and $\widetilde{D}_{h,t}$ is the average net demand for each home $h$ at time $t$ obtained from the historical data (consumption and PV generation) analysis [31]; $\widetilde{Loss}_t$ is the total active losses of the network lines at time $t$ obtained from the load flow program; $P_t^{DG}$ is the available power from the distributed energy resources owned by the DNO.

For the day-ahead cost function $F_t^{DHA}$, the well-known quadratic model which is a non-decreasing convex function in $g_t$ is used [32] as in (3).

$$F_t^{DHA}(g_t) = a g_t^2 + b g_t + c \tag{3}$$

where $a$,$b$ and $c$ are constants.

Please note the proposed framework is general and any method can be used for the day-ahead price calculations as well as the day-ahead cost function modelling within this framework.

## 4.2    Model of customers

The HEMS controls the operation of home devices aiming to find the best operation times for appliances and charging/discharging profile of the BESS. It starts to operate once the day-ahead price signals or consumption adjustment requests are received. Then, the total hourly consumptions of the home are sent back to the DNO. The HEMS model developed in [31] is considered here, which is explained briefly in Appendix A.

The total electricity consumption of home $h$ at time $t$ is represented by $L_{h,t}$ which is calculated based on the operational status of home devices, PV generation and BESS power. The positive/negative $L_{h,t}$ corresponds to the consumption/injection from/into the network. The goal of HEMS is to decide when to turn on/off home devices and find the best operational mode of the battery over the day to minimize the customer's electricity payment. Each HEMS solves an optimization problem with the objective function $J_h(.)$ given by (4):

$$J_h(L_{h,t}) = \min_{L_{h,t}} \sum_{t=1}^{T} \lambda_t L_{h,t} \Delta t \tag{4}$$

where

$$\lambda_t = \begin{cases} \pi_t, & L_{h,t} > 0 \\ FiT_t, & L_{h,t} < 0 \end{cases} \tag{5}$$

Where $\Delta t$ is the duration of each time slot and $T$ is the total number of time slots. $\pi_t$ and $FiT_t$ are the electricity price and FiT rate, respectively.

## 4.3    Distribution network operator model

The operator is responsible for congestion management by designing electricity prices in a way that encourages customers to schedule their consumption to alleviate congestion in the network. The profit value for the operator at each time $t$ is approximately computed by:

$$Profit_t = \sum_{h=1}^{H} J_h(L_{h,t}) - F_t^{DHA} - F_t^{RT} \tag{6}$$

The first term is the sum of customer's electricity payments which represents the income of the DNO. Please note that the payment of DNO for the customers feed-in-electricity is already included in $\sum_{h=1}^{H} J_h(L_{h,t})$. $F_t^{DHA}$ represents the day-ahead market cost defined in (3) and $F_t^{RT}$ models the cost function of the operator in the real-time balance market described by:

$$F_t^{RT} = \mu_t^b \left( \sum_{h=1}^{H} L_{h,t} + Loss_t - g_t \right)^+ - \mu_t^s \left( \sum_{h=1}^{H} L_{h,t} + Loss_t - g_t \right)^- \tag{7}$$

where $g_t$ is the total forecasted demand of the network (i.e., the power already purchased from the day-ahead market); $\mu_t^b$ and $\mu_t^s$ denote the price of purchasing and absorbing extra power in the real-time market, respectively. $Loss_t$ is the total active losses of the network in real-time operation.

The individual operation of multiple HEMS will possibly increase $F_t^{RT}$. Thus, it is profitable for the operator to focus on keeping the actual load as close to the available supply as possible. This can also potentially benefit customers in a long-term view since more profit for the operator means that it might reduce the electricity prices. Therefore, minimization of $F_t^{RT}$ is considered as the objective of the proposed coordinated load scheduling model at the network level.

### 4.4 Centralized model for the coordinated load scheduling

Considering the desires of all players, the overall objective function of the coordinated optimization model is:

$$\min_{L_{h,t}} \sum_{t=1}^{T} F_t^{RT}(\boldsymbol{L}_t^{Toal}) + \alpha \sum_{t=1}^{T} \sum_{h=1}^{H} J_h(L_{h,t}) \tag{8}$$

*Subject to:*

$$\boldsymbol{L}_t^{Toal} = \sum_h L_{h,t} + Loss_t, \qquad \forall t \in T \tag{9}$$

$$I_{f,min} \leq I_{f,t} \leq I_{f,max}, \qquad \forall t \in T, \forall f \in N_F \tag{10}$$

$$L_{h,t} \in \mathcal{M}_{h,t}, \qquad \forall h \in H, \forall t \in T \tag{11}$$

where $\alpha$ is the scaling factor; $I_{f,t}$ the current of line $f$ at time $t$; and $N_F$ is the set of all network lines. The first term in (8) represents the cost of DNO from participating in the real-time balancing market and the second term is the sum of individual payments of customers. The coordinated optimization model is subject to both network-related and home-related constraints. Equation (9) determines the energy balance constraint within the network including network losses which is obtained by adding active losses of all lines. Line current constraints are included in (10) obtained from the load flow program using the forward-backward sweep method [33]. Constraint (11) and $\mathcal{M}_{h,t}$ characterizes all the operational constraints of devices inside the home which are defined in detail in [31]. They include desired operating period, length of operation and power consumption of appliances and (dis)charging limits of the BESS. In the proposed model, customers do not experience any degree of discomfort because their preference in running appliances is included in the formulation through (11). The coordinated model only shifts the activation time of appliances over the day for a better load profile enhancement, but the operation time is still within the desired period of the users.

## 4.5 Distributed optimization for the coordinated load scheduling

The formulated optimization model in (8)-(10) gives the optimal solution of the problem, however, it is not practically feasible to solve it in a centralized manner due mainly to the huge dimensionality and privacy-related issues. Specifically, with the increase in the number of customers, the number of variables and constraints will prohibitively enlarge which make solving the optimization problem insufficient. It has been shown that the average runtime for a centralized optimal power flow problem will grow exponentially or linearly with a high slope. Moreover, most of the existing hardware setups may fail to solve such a centralized model because of the memory overflow [34, 35].

A distributed optimization method is preferred, but (8) is not solvable in a distributed manner in the current form due to the existence of the coupling variable $L_{h,t}$ between the DNO and customers. This is a sharing problem that can be resolved using the well-known ADMM technique, by decomposing the coupling constraints (Appendix B). The main idea behind the distributed load scheduling is to enable customers to autonomously optimize their consumption in response to the DNO signals. In this section, the centralized coordination model is modified to design a distributed load scheduling model.

To be able to solve the centralized problem by ADMM algorithm, a set of auxiliary variables $\hat{L}_{h,t}$ for the total electricity consumption of home $h$ at time $t$ is introduced. Then, the problem can be formulated as:

$$\min_{L_{h,t}, \hat{L}_{h,t}} \sum_{t=1}^{T} F_t^{RT}(\hat{L}_t^{Toal}) + \alpha \sum_{t=1}^{T} \sum_{h=1}^{H} J_h(L_{h,t}) \tag{12}$$

$subject\ to$:

$$\hat{L}_t^{Toal} = \sum_h \hat{L}_{h,t} + \widehat{Loss}_t, \qquad \forall t \in T \tag{13}$$

$$I_{f,min} \leq I_{f,t} \leq I_{f,max}, \qquad \forall t \in T, \forall f \in N_F \tag{14}$$

$$L_{h,t} \in \mathcal{M}_{h,t}, \qquad \forall h \in H, \forall t \in T \tag{15}$$

$$L_{h,t} = \hat{L}_{h,t}, \qquad \forall t \in T, \forall h \in H \tag{16}$$

The first term is the objective function of the DNO and the second term includes the objective function of the customers. Each term has separate decision variables and subjected to separate constraints, but the constraint $L_{h,t} = \hat{L}_{h,t}$ links them together. The distributed version of the problem includes two groups of optimization variables; $L_{h,t}$ is the perception of customer $h$ of its consumption at time $t$ while $\hat{L}_{h,t}$ captures the DNO perception of the consumption of customer $h$ at time $t$. For the DNO, constraints are included in (13)-(14), while each customer $h$ aims to satisfy their constraints in (15). Through the distributed optimization process, variables of each group are optimized individually in a parallel manner and are coordinated together via the associate dual variables $u_{h,t}$. Ultimately, $L_{h,t} = \hat{L}_{h,t}$ is achieved.

By applying the ADMM algorithm to the problem in (12)-(16), the augmented Lagrangian function is calculated using (17):

$$\begin{aligned} &\mathcal{L}_\rho(\{L_{h,t}\}, \{\hat{L}_{h,t}\}, \{u_{h,t}\}) \\ &= \sum_{t \in T} F_t^{RT}(\hat{L}_t^{Toal}) + \alpha \sum_{t \in T} \sum_{h \in H} J_h(L_{h,t}) + \frac{\rho}{2} \sum_{t \in T} \sum_{h \in H} \left\| L_{h,t} - \hat{L}_{h,t} + u_{h,t} \right\|_2^2 \end{aligned} \tag{17}$$

where $u_{h,t}$ is the scaled dual variable associated with home $h$ at time $t$.

Please note that only the coupling constraint $L_{h,t} = \hat{L}_{h,t}$ is included in the Lagrangian formulation. The DNO and each customer optimize their objective functions and update their decision variables. Thereby, at each iteration $k + 1$, the following steps are done by each player:

| Algorithm 1: Distributed algorithm for a coordinated load scheduling |
|---|
| 1: Inputs: $\mu_t^p, \mu_t^s, \lambda_t, \{P_{h,t}^{pv}\}_{t=1,h=1}^{T,H}, \{P_{h,t}^{non}\}_{t=1,h=1}^{T,H}, \boldsymbol{\Phi}_{h,t}^a, \boldsymbol{\Gamma}_{h,t}^a, \rho, \alpha$ |
| 2: Initialize: $\{L_{h,t}^0\}_{t=1,h=1}^{T,H}, \{\hat{L}_{h,t}^0\}_{t=1,h=1}^{T,H}, \{u_{h,t}^0\}_{t=1,h=1}^{T,H}, k = 0$; |
| 3: **while** not converged **do** |
| 4:    $k = k + 1$; |
|      # <u>at the customer side</u> # |
| 5:    **for** $h = 1:H$ and in parallel for $h$ **do** |
| 6:      Solve (18) and obtain $\{L_{h,t}^k\}_{t=1}^T$ and send them to the DNO |
| 7:    **end for** |
|      # <u>at the DNO side</u> # |
| 8:    Wait to receive $\{L_{h,t}^k\}_{t=1,h=1}^{T,H}$ from all customers |
| 9:    Solve (19) and find $\{\hat{L}_{h,t}^k\}_{t=1,h=1}^{T,H}$ for each customer |
| 10: Update dual variables by: $u_{h,t}^k = u_{h,t}^{k-1} + L_{h,t}^k - \hat{L}_{h,t}^k, \forall h \in H, \forall t \in T$ |
| 11: Send $\{\hat{L}_{h,t}^k\}_{t=1,h=1}^{T,H}$ and $\{u_{h,t}^k\}_{t=1,h=1}^{T,H}$ to each customer |
| 12: **end while** |
| 13: Outputs: $\{L_{h,t}^k\}_{t=1,h=1}^{T,H}$ |

- Customer $h$ updates:

$$\{L_{h,t}^{k+1}\}_{t=1}^T = \underset{\{L_{h,t}\}_{t=1}^T}{\text{argmin}} \ \alpha \sum_{t \in T} J_h(L_{h,t}) + \frac{\rho}{2} \sum_{t \in T} \left\| L_{h,t} - \hat{L}_{h,t}^k + u_{h,t}^k \right\|_2^2 \tag{18}$$

- DNO update:

$$\{\hat{L}_{h,t}^{k+1}\}_{t=1,h=1}^{T,H} = \underset{\{\hat{L}_{h,t}^{k+1}\}_{t=1,h=1}^{T,H}}{\text{argmin}} \ \sum_{t \in T} F_t^{RT}\left(\sum_{h \in H} \hat{L}_{h,t}\right) + \frac{\rho}{2} \sum_{h \in H} \sum_{t \in T} \left\| L_{h,t}^{k+1} - \hat{L}_{h,t} + u_{h,t}^k \right\|_2^2 \tag{19}$$

- Dual variables updates:

$$u_{h,t}^{k+1} = u_{h,t}^k + L_{h,t}^{k+1} - \hat{L}_{h,t}^{k+1}, \forall h \in H, \forall t \in T \tag{20}$$

In (18), $\{L_{h,t}\}_{t=1}^T$ are updated by solving a load scheduling problem for each home while considering $\hat{L}_{h,t}$ as a fixed value from the previous iteration. This optimization problem is subjected to the constraint (15) which characterizes all the operational constraints of devices inside home $h$. Therefore, the customer's comfort level in terms of appliance operation time is preserved. Customers solve their optimization problems and update their consumption in a parallel fashion. Then, they only need to transfer their overall consumption at each time, i.e., $L_{h,t}$, to the DNO. Intuitively, $L_{h,t}$ is the load of customer $h$ for her benefit. The DNO then solves the optimization problem in (19) and finds $\hat{L}_{h,t}$ for each home $h$ at each time $t$ taking network-related constraints in (13)-(14) into account. Intuitively speaking, $\hat{L}_{h,t}$ can be interpreted as the suggested load by the DNO for home $h$ at time $t$ for the safe operation of the network and the DNO's cost minimization. Finally, $L_{h,t}$ and $\hat{L}_{h,t}$ are coupled together by the scaled dual variable in (20). The DNO updates $u_{h,t}$ for the consumption of each home at each time slot and announces to each customer for the next iteration. Putting it differently, the whole algorithm is the process of negotiation between the DNO and each customer.

A step-by-step explanation of the distributed load scheduling is shown in Algorithm 1. The centralized large-scale problem is now decomposed into $H + 1$ small-scale problems; one corresponding to the DNO and the rest associated with each customer. There is also a dual variable update phase which is done by the DNO. The proposed distributed energy management starts by announcing the day-ahead price values by DNO to the customers' HEMS. Then, each HEMS finds the best schedule of appliances individually and sends its overall consumption at each time of the day to the DNO. Then, once all customers send their consumption information, the DNO runs the load flow program and finds new values for the consumption of each home intending to minimize its cost in the real-time market and maintaining network's flows inside the operating limits. Then, dual variables are updated along with the desired consumption values are sent back to each customer where each HEMS updates the schedule of appliances to meet the needs of the network as much as possible. Dual variables are coordinating the load schedule between the DNO and customers. This process is repeated until a convergence happens and the DNO and customers reach an agreement on the electricity consumption for each time of the day. The interaction between customers and the DNO in reaching an agreement on the consumption at each time $t$ is depicted in Figure 2.

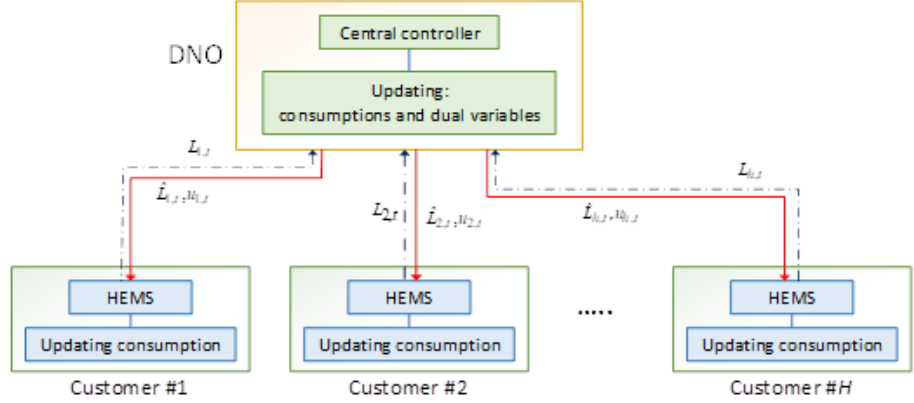

Figure 2: Interaction between the DNO and customers in the proposed coordinated load scheduling model

### 4.6 Incentive design

The price adjustment policy is to guide the electricity consumption behavior of the customers by raising the price for eliminating unwanted peak loads and reducing the price for filling the valley load. In other words, customers get rewarded for decreasing their consumption on peak loads or increasing it to fill the valley loads. Similarly, they get penalized for increasing their consumption on peak loads or decreasing it on valley loads. The price adjustment policy is explained in Figure 3. In this regard, time-dependent rewards and penalties are added to the base price value and offered to all customers intending to incentivize them to follow the desired consumption of the DNO. To this end, a new price term $\in_t$ is calculated and added to the base price $\lambda_t$ which depends on the difference between the actual load of the network and one that is requested by the DNO and is desirable for the safe operation of the network.

The new price term $\in_t$ is a function of the cost that the DNO incurs in the real-time balancing market. In this way, the retail electricity price changes as the real-time wholesale price changes which guarantee the benefit for the operator. Therefore, the price gap $\in_t$ is calculated as in (21):

$$\in_t = \sigma_t \frac{\partial F_t^{RT}(\Delta L_t)}{\partial \Delta L_t} \qquad (21)$$

where $\sigma_t$ is the reward/penalty coefficient that captures the influence of the deviation between the actual and desired load. $\Delta L_t$ denotes the difference between the actual and desired load within the network at time $t$ which is obtained by:

$$\Delta L_t = \hat{\boldsymbol{P}}_t - \boldsymbol{P}_t^{unc} \qquad (22)$$

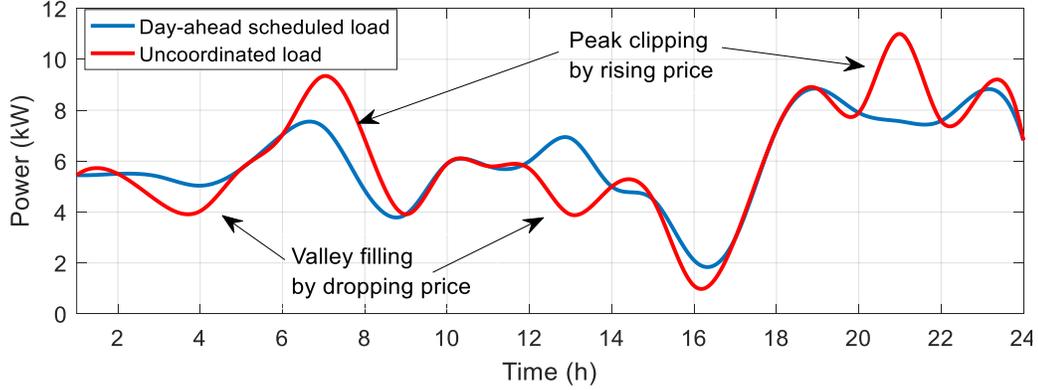

Figure 3: Price adjustment policy in the coordinated load scheduling program

where $\hat{P}_t = \sum_{h \in H} \hat{L}_{h,t} + \widehat{Loss}_t$ is the total load of the network preferred by DNO at time $t$ and $\widehat{Loss}_t$ is the associated network losses. Ideally suited, $\hat{P}_t = g_t$. $P_t^{unc} = \sum_{h \in H} L_{h,t}^{unc} + Loss_t^{unc}$ is the total load for the uncoordinated and individual operation of HEMS based on the base price values, where $L_{h,t}^{unc}$ and $Loss_t^{unc}$ are the original load of customer $h$ and network losses at each time $t$, respectively.

The new price term $\in_t$ should also be a function of the customer's contribution to load adjustment. To achieve this, the reward/penalty coefficient $\sigma_t$ is calculated based on the total shifted energy as [3]:

$$\sigma_t = \begin{cases} w^{-1}\left(e^{\left(\frac{1}{w}\right)\left(\frac{\hat{P}_t - P_t^{unc}}{P_t^{unc}}\right)} - 1\right), & \Delta L_t > 0 \\ w\left(e^{\left(\frac{1}{w}\right)\left(\frac{\hat{P}_t - P_t^{unc}}{P_t^{unc}}\right)} - 1\right), & \Delta L_t < 0 \end{cases} \tag{23}$$

where $w \in [0,1]$ is a parameter that tunes the amount of reward or penalty added to the base price. It controls how generous and strict the DNO is in rewarding or penalizing customers. Its influence will be demonstrated later by numerical examples.

The reward/penalty coefficient $\sigma_t$ is modelled using an exponential function to provide increasing reward or penalty for the participation of customers and the value of $\sigma_t$ is non-zero only in the time slots when there is a deviation between the aggregated load in response to the day-ahead price, i.e., $P_t^{unc}$, and the desired load for the network operation, i.e, $\hat{P}_t$. In this manner, when $\hat{P}_t > P_t^{unc}$ and accordingly $\Delta L_t > 0$, then $\sigma_t$ becomes negative which leads to a reduction in the price motivating customers to increase their consumption. Conversely, when $\hat{P}_t < P_t^{unc}$ and accordingly $\Delta L_t < 0$, which indicates that the aggregated load is more than the available supply. Therefore, the value of $\sigma_t$ becomes positive, and consequently, the price is increased to penalize customers forcing them to reduce their consumption.

The reward/penalty coefficient $\sigma_t$ is a function of the overall contribution of customers in reducing the mismatch between the available supply and aggregated load, which motivates customers to schedule their loads cooperatively to minimize this mismatch as much as possible. It is also worth mentioning that the new price value, $\in_t$, is added to the price of selling electricity to the grid $Fit_t$. In this case, rewarding customers means to increase the FiT rate whereas reducing it corresponds to penalizing them. In this regard, when available supply is higher than the demand, i.e., $\hat{P}_t > P_t^{unc}$, $\sigma_t$ becomes negative and hence FiT rate is reduced to prevents customers from injecting power into the grid. Likewise, when demand is higher than the available supply, i.e., $\hat{P}_t < P_t^{unc}$, $\sigma_t$ becomes positive leading to an increase in the FiT rate. This motivates generating customers to inject their surplus power into the grid to make more money for themselves and at the same time help the DNO in maintaining a power balance in the network. However, considering that in the proposed load scheduling model, more attention is given to rewarding customers rather than penalizing them, and to keep financially incentivizing residential customers to install solar PV systems, the penalty value related to the FiT price is zero. In other words, the FiT price will never fall below its previously announced value, but it might be increased as an incentive for more power injection. Given this, the reward/penalty coefficient $\sigma_t$ related to the FiT price is calculated as in (24).

$$\sigma_t = \begin{cases} 0 \,, & \Delta L_t > 0 \\ w^{-1}\left(e^{\left(\frac{1}{w}\right)\left(\frac{\bar{P}_t - P_t^{unc}}{P_t^{unc}}\right)} - 1\right), & \Delta L_t < 0 \end{cases} \quad (24)$$

Therefore, the real-time price offered to customers in the proposed coordinated load scheduling model is obtained as in (25).

$$\lambda_t^{RT} = \lambda_t + \epsilon_t \quad (25)$$

Then, by substituting (25) in (8) and computing the payment of customers for the real-time price, the objective function of the centralized model is calculated using (26).

$$\min_{L_{h,t}} \sum_{t=1}^{T} F_t^{RT}(\boldsymbol{L}_t^{Toal}) + \alpha \left[\sum_{t=1}^{T}\sum_{h=1}^{H} J_h(L_{h,t}) + \sum_{t=1}^{T}\sum_{h=1}^{H} \epsilon_t\, L_{h,t}\right] \quad (26)$$
$$subject\ to: \qquad (13) - (16)$$

By employing the auxiliary variables $\hat{L}_{h,t}$, the equivalent objective function is obtained as in (27).

$$\min_{L_{h,t},\hat{L}_{h,t}} \sum_{t=1}^{T}\left[F_t^{RT}(\hat{\boldsymbol{L}}_t^{Toal}) + \alpha \sum_{h=1}^{H} \epsilon_t\, \hat{L}_{h,t}\right] + \alpha \sum_{t=1}^{T}\sum_{h=1}^{H} J_h(L_{h,t}) \quad (27)$$
$$subject\ to: \qquad (13) - (16)$$

By applying the ADMM algorithm, the Lagrangian function is calculated using (28).

$$\mathcal{L}_\rho(\{L_{h,t}\}, \{\hat{L}_{h,t}\}, \{u_{h,t}\})$$
$$= \sum_{t\in T}\left[F_t^{RT}(\hat{\boldsymbol{L}}_t^{Toal}) + \alpha \sum_{h=1}^{H} \epsilon_t\, \hat{L}_{h,t}\right] + \alpha \sum_{t\in T}\sum_{h\in H} J_h(L_{h,t})$$
$$+ \frac{\rho}{2}\sum_{t\in T}\sum_{h\in H}\|L_{h,t} - \hat{L}_{h,t} + u_{h,t}\|_2^2 \quad (28)$$

Note that the newly added term is only included in the DNO objective function. Therefore, it results in a different equation for the DNO in the ADMM-based distributed algorithm. Thus, at iteration $k + 1$, the DNO solves the following optimization problem:

- DNO update:

$$\{\hat{L}_{h,t}^{k+1}\}_{t=1,h=1}^{T,H} = \underset{\{\hat{L}_{h,t}^{k+1}\}_{t=1,h=1}^{T,H}}{\operatorname{argmin}} \sum_{t\in T}\left[F_t^{RT}\left(\sum_{h\in H}\hat{L}_{h,t}\right) + \alpha \sum_{h\in H} \epsilon_t\, \hat{L}_{h,t}\right]$$
$$+ \frac{\rho}{2}\sum_{h\in H}\sum_{t\in T}\|L_{h,t}^{k+1} - \hat{L}_{h,t} + u_{h,t}^k\|_2^2 \quad (29)$$

The new term added to the DNO problem in (29) is associated with the reward/penalty amount for each customer at each hour of the day for their participation in the DR program. At each iteration of the algorithm, the DNO computes the reward/penalty values for customers based on their response to the requested power in the previous iteration and tries to minimize the total sum of rewards or penalties of all customers for the whole scheduling period (i.e., a day). At the end of the scheduling period, each

---

**Algorithm 2: Distributed algorithm for a coordinated load scheduling with price adjustment included**

---

1: Inputs: $\mu_t^b, \mu_t^s, \lambda_t, \{P_{h,t}^{pv}\}_{t=1,h=1}^{T,H}, \{P_{h,t}^{non}\}_{t=1,h=1}^{T,H}, \boldsymbol{\Phi}_{h,t}^a, \boldsymbol{\Gamma}_{h,t}^a, \rho, \alpha$

2: Initialize: $\{L_{h,t}^0\}_{t=1,h=1}^{T,H}, \{\hat{L}_{h,t}^0\}_{t=1,h=1}^{T,H}, \{u_{h,t}^0\}_{t=1,h=1}^{T,H}, k = 0;$

3: **while** not converged **do**

4: $\quad k = k + 1;$

*\* at the customer side \**

5: $\quad$ **for** $h = 1:H$ and in parallel for $h$ **do**

6: $\quad\quad \underset{\{L_{h,t}\}_{t=1}}{\text{minimize}} \quad \alpha \sum_{t \in T} J_h(L_{h,t}) + \frac{\rho}{2} \sum_{t \in T} \left\| L_{h,t} - \hat{L}_{h,t}^k + u_{h,t}^k \right\|_2^2$

$\quad\quad\quad\quad\quad\quad\quad\quad\quad\quad\quad\quad\quad subject\ to:\ (15)$

7: $\quad\quad$ send $\{L_{h,t}^k\}_{t=1}^{T},$ to the DNO

8: $\quad$ **end for**

*\* at the DNO side\**

9: $\quad$ Wait to receive $\{L_{h,t}^k\}_{t=1,h=1}^{T,H}$ from all customers

10: $\quad \underset{\{L_{h,t}^{k+1}\}_{t=1,h=1}^{T,H}}{\text{minimize}} \quad \sum_{t \in T}\left[F_t^{RT}\left(\sum_{h \in H} \hat{L}_{h,t}\right) + \alpha \sum_{h \in H} \in_t \hat{L}_{h,t}\right] + \frac{\rho}{2} \sum_{h \in H} \sum_{t \in T} \left\| L_{h,t}^{k+1} - \hat{L}_{h,t} + u_{h,t}^k \right\|_2^2$

$\quad\quad\quad\quad\quad\quad\quad\quad\quad\quad\quad subject\ to\ (13) - (14)\ \&\ (31)$

11: $\quad$ Update dual variables by: $u_{h,t}^k = u_{h,t}^{k-1} + L_{h,t}^k - \hat{L}_{h,t}^k, \forall h \in H, \forall t \in T$

12: $\quad$ Send $\{\hat{L}_{h,t}^k\}_{t=1,h=1}^{T,H}$ and $\{u_{h,t}^k\}_{t=1,h=1}^{T,H}$ to each customer

13: **end while**

14: Outputs: $\{\in_t\}_{t=1}^T,$ and $\{L_{h,t}\}_{t=1,h=1}^{T,H}$

---

customer will receive its rebate on the electricity bill in the form of adjustments on the electricity prices. She might get rewarded at some hours and might get penalized at some other hours during the day.

The incentives given to the customers should be non-discriminatory. Intuitively, this means that they all get the same amount of rebate because of their contribution to shaping the network's load profile. Hence, a constraint is added to the optimization problem solved by the DNO. In this regard, the normalized rebate $R_h$ given to each customer $h \in H$ is calculated using (30).

$$R_h = \frac{J_h\left(L_{h,t}^{unc}\right) - J_h\left(\hat{L}_{h,t}\right)}{J_h\left(L_{h,t}^{unc}\right)} \tag{30}$$

where $J_h\left(L_{h,t}^{unc}\right)$ is the daily payment of customer for the uncoordinated and individual operation based on the base price values; $J_h\left(\hat{L}_{h,t}\right)$ is the daily payment of customer for following the suggested load by the DNO, i.e., $\hat{L}_{h,t}$, based on the adjusted electricity price values. Then, the constraint (31) is added.

$$\sum_{h=1}^{H} \left| R_h - R_{avg} \right| \leq \varepsilon \tag{31}$$

where $R_{avg}$ is the average value of all customer's rebates. Ideally, customers should be receiving the same amount of reduction or decrease in their daily payments and thus $\varepsilon = 0$. However, because customers have different daily consumption patterns it is not practically possible to reach zero for $\varepsilon$, but it is normally chosen to be a small value. This constraint is included in the model in (29) solved by the DNO.

A detailed description of the proposed ADMM-based load scheduling algorithm highlighting the customer and DNO optimization steps is summarised in Algorithm 2. The optimization problem of the DNO is a non-linear programming (NLP) model which is solved in MATLAB using "fmincon" solver. For the customer optimization problem which is a mixed-integer non-linear programming (MINLP) model, the proposed HEMS model in [31] has been used.

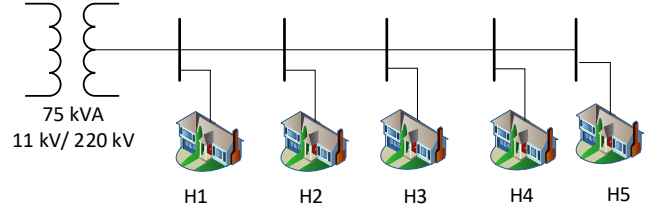

Figure 4: Schematic of the residential test feeder

Table 1: Desired time interval of each customer for the operation of appliances

| | Controllable appliance number | | | | | |
|---|---|---|---|---|---|---|
| | #1 | #2 | #3 | #4 | #5 | #6 |
| Home #1 Home #2 Home #5 | 4 pm- 11pm | 10 am- 4 pm | 5 am- 9 am | 5 am- 9 am | 10 am- 4 pm | 4 pm- 11 pm |
| Home #4 Home #5 | 9 am- 8 pm | 9 am- 8 pm | 9 am- 8 pm | 9 am- 8 pm | 9 am- 8 pm | 9 am- 8 pm |

## 5 Numerical Simulations

### 5.1 Data preparation and parameters setting

The proposed DR model is implemented in a simple test residential feeder with five customers ( Figure 4) and associated data in Appendix C). These homes are selected from the "Ausgrid" database [36]. Among them, home #3 and home #4 have installed solar PV with a size of 6 and 4.7 kWp, respectively. Also, home #3 has installed a BESS with a size of 5 kWh. It is assumed that each home has 6 controllable appliances each of which with different consumption profiles and 1-hour job length to finish the task for simplicity. The energy usage of these appliances is assumed to be 0.4, 0.8, 1.2, 1.6 2, and 2.5 kWh. Each customer might have a different desired time for activating each of these appliances which are shown in Table 1. The scheduling period is set to a day ($T = 24$) with a one-hour time resolution ($\Delta t = 1\ hr$).

To obtain a forecast of the network demand for the next day, the average net consumption for homes at each hour of the day is found using data analysis [31]. It is also assumed that the DNO does not own any energy resources. Thus, the forecasted demand for the day and the day-ahead base prices are calculated and shown in Figure 5. For price gap value $\sigma_t$ in (23) and (24), $w = 0.5$ is chosen. Other simulation parameters are listed in Appendix C.

The convergence behavior of the ADMM algorithm is very sensitive to the value of the penalty parameter $\rho$. A very popular technique is to change $\rho$ iteratively based on the performance of the algorithm as follows [28, 37]:

$$\rho^{k+1} = \begin{cases} \tau^{incr}\rho^k, & if\ \|r^k\|_2 > \gamma\|s^k\|_2 \\ \rho^k/\tau^{decr} & if\ \|s^k\|_2 > \gamma\|r^k\|_2 \\ \rho^k, & otherwise, \end{cases} \tag{32}$$

where $\gamma > 1$, $\tau^{incr} > 1$, and $\tau^{decr} > 1$ are parameters, $r^k$ is the primal residual, and $s^k$ the dual residual of the algorithm. The main goal is to vary $\rho$ in line with the relative magnitude of the primal and dual residential to make $\rho^k \to \infty$ and accordingly, keep the primal and dual residential norms within a factor of $\gamma$ of each other as they converge to zero. In this paper, the suggested values in [28] as $\gamma = 10$, $\tau^{incr} = \tau^{decr} = 2$ and $\rho^0 = 0.001$ are used. Also, the value of the scaled dual variables $u$ is rescaled after updating $\rho$. This formula improves the convergence and makes algorithm performance less sensitive to the initial value of $\rho$ [37, 38].

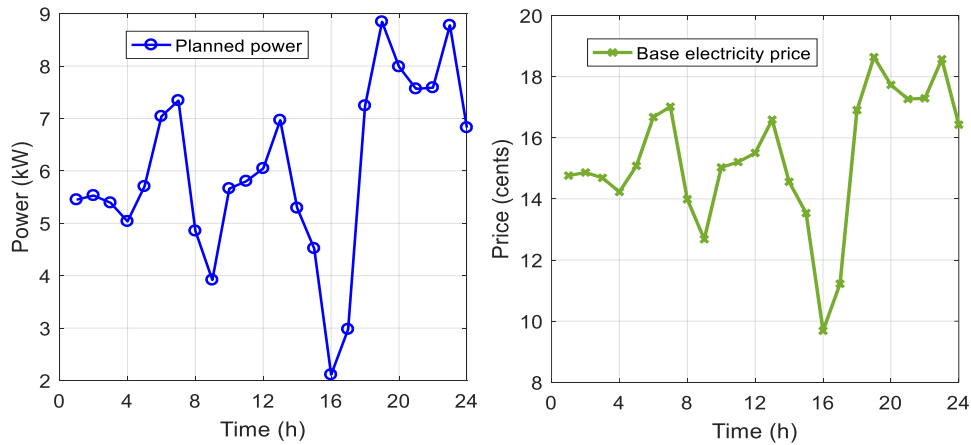

Figure 5: Day-ahead planned power and designed price values for the network by DNO

### 5.2    Performance evaluation

The aggregated load of the network with uncoordinated and coordinated HEMS operation (proposed in this paper) against the planned load for the network in a given day is shown in Figure 6.  When customers use their HEMS individually without any coordination, a significant power imbalance and a huge deviation from the day-ahead planned load are observed. Specifically, as expected, the two rebound peaks during the low-pricing hours at 8 a.m. and 3-4 p.m. are observed, which results in a high real-time cost for the DNO and possible violations in the network constraints. The coordinated HEMS approach as proposed in this paper, however, helps in alleviating the rebound peaks resulting in a flatter load profile with less deviation from the day-ahead load. This is achieved by adjusting the day-ahead price values to incentivize customers to shift their consumption from low-pricing hours to other hours of the day.

The adjusted price values are depicted in Figure 7. The price of electricity is increased at 8 a.m., 1 p.m., 3 p.m., and 4 p.m., i.e., whenever the actual consumption is higher than the available supply to penalize customers. Simultaneously, the FiT price is increased as a reward for generating customers to inject more power into the gird minimizing net power imbalances. On the other hand, the price of electricity is reduced at times when the consumption is less than the available supply to reward customers for shifting their consumption to these hours. It is ensured that the customers are not penalized for injecting more power into the grid and thus the FiT rate is not reduced.

As a result of coordinating HEMS operation, unwanted peaks eliminated, and the valley loads filled and accordingly the cost of the DNO in the real-time market has dramatically fallen from 118.11 to 30.96 cents (75% saving). Customer's electricity bills have decreased too. The electricity payments of the customers before and after participating in the proposed load scheduling model are compared in Table 2. All customers received around 4.5% reduction on their daily electricity bill. Therefore, the proposed coordinated operation benefits both customers and the DNO.

Furthermore, the proposed approach reduces total network losses, by 3.5% (from 3.129 kW to 3.022 kW). To further quantify the improvement in the network load profile, the peak-to-average (PAR) ratio for the aggregated load of the network is calculated by dividing the network's peak load by the average load as in (33) [39].

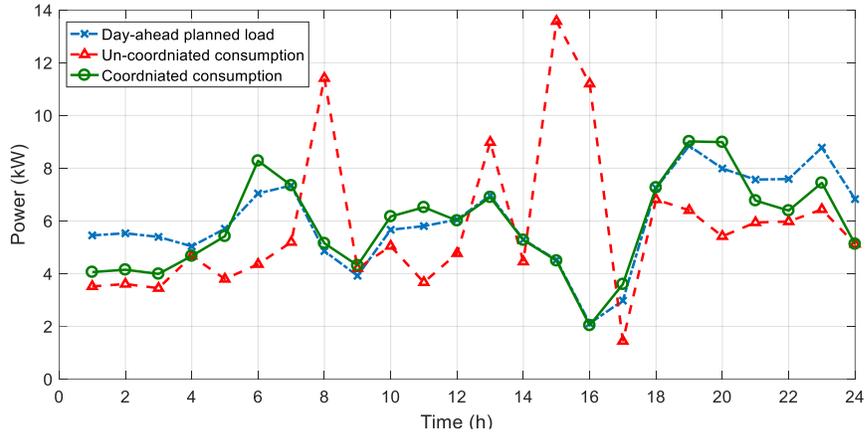

Figure 6: Comparison of network load before and after implementing the proposed model

$$PAR = \frac{\boldsymbol{P}_{peak}}{\boldsymbol{P}_{avg}} = \frac{T \max_{t \in T} \boldsymbol{P}_t}{\sum_{t \in T} \boldsymbol{P}_t} \qquad (33)$$

The PAR of the network has decreased from 2.336 for the uncoordinated operation to 1.552 for the coordinated operation, which shows around 34% enhancement.

Figure 8 demonstrates the benefits of mitigating the loading of the network lines. With the uncoordinated and individual operation of HEMS, high line overloading is observed resulting in an increased damage cost for the DNO. By coordinating the operation of HEMS, overloading is eliminated resulting in potential savings for DNO.

Figure 9 displays the detailed load consumption of home #3 and #4. The appliances in each home were scheduled to operate during the sun-shining and low-price hours, i.e., 9 am to 4 pm. This is intuitive for individual HEMS as they want to maximize the usage of in-house PV generation and low-price hours. Whereas, once the proposed coordinated HEMS operation is used, individual HEMS are requested to inject their PV generation into the grid at 1-4 pm as much as possible to help reduce the network load. In return, the FiT price is increased to reward them. The BESS profile of home #3 for two cases is also compared in Figure 10. During the uncoordinated operation, the battery is charged by the cheapest electricity in the morning to get discharged later for covering the household's morning peak at 6 am. Also, it is fully charged from the PV power during the day to supply household demand later in the evening. The battery profile, however, has slightly changed because of the participation of the customer in the DR program. The battery is now charged between 1-4 am when the total demand of the network is less than the available supply. The price of electricity is decreased at these times to benefit the customer's participation. The customer is also required to discharge its battery at 4 pm when the available supply is much less than the total demand. The customer, in exchange, is paid back on a higher FiT rate.

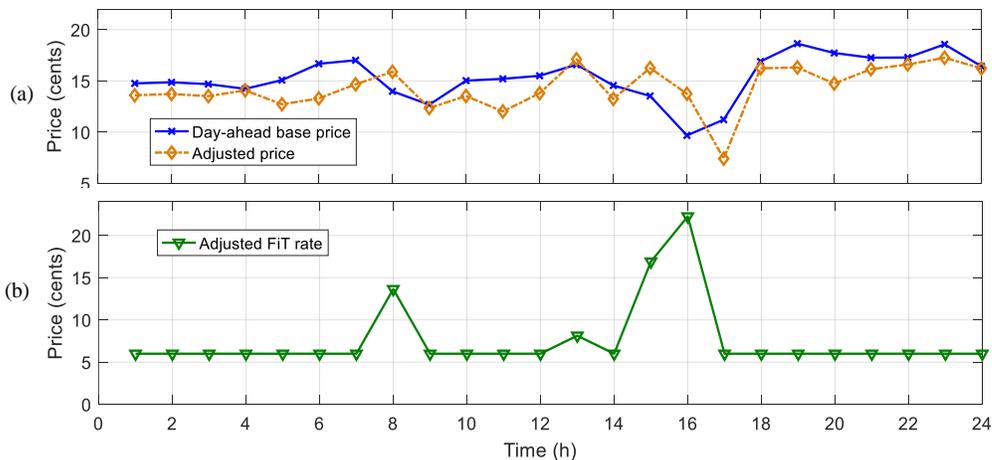

Figure 7: Adjusted price with reward/penalty values, (a) price of buying electricity from the grid, (b) FiT price

Table 2: Payment comparison of customers for their participation in the proposed load scheduling model

| | Electricity payment (cents) | | Saving |
|---|---|---|---|
| | uncoordinated operation | coordinated operation | |
| Home #1 | 405.872 | 387.358 | 4.56 % |
| Home #2 | 792.959 | 753.714 | 4.94 % |
| Home #3 | 149.847 | 143.598 | 4.17 % |
| Home #4 | 357.306 | 340.951 | 4.57 % |
| Home #5 | 475.961 | 459.423 | 3.47 % |

### 5.3 Convergence behavior

Figure 11 demonstrates the load modification progress of customers at some hours of the day in the coordinated load scheduling model. It shows how customers and the DNO iteratively negotiate for the hourly consumption and finally reach an agreed value. The DNO and customers begin with different consumption values, but as the coordination process starts and they communicate their desired values, they move toward each other in small steps until they come to an agreement and achieve a common value. The agreed consumption value is desirable for the DNO and preserves the customer's comfort.

Figure 12 shows the change in the total electricity bill of customers versus the iteration number. The evolution of primal and dual residuals of the ADMM algorithm are also displayed. The algorithm reaches a convergence after around 50 iterations. It is also worth noting that the distributed optimization is done in a parallel fashion meaning that each customer sends/receives consumption values to/from the DNO at the same time with other customers.

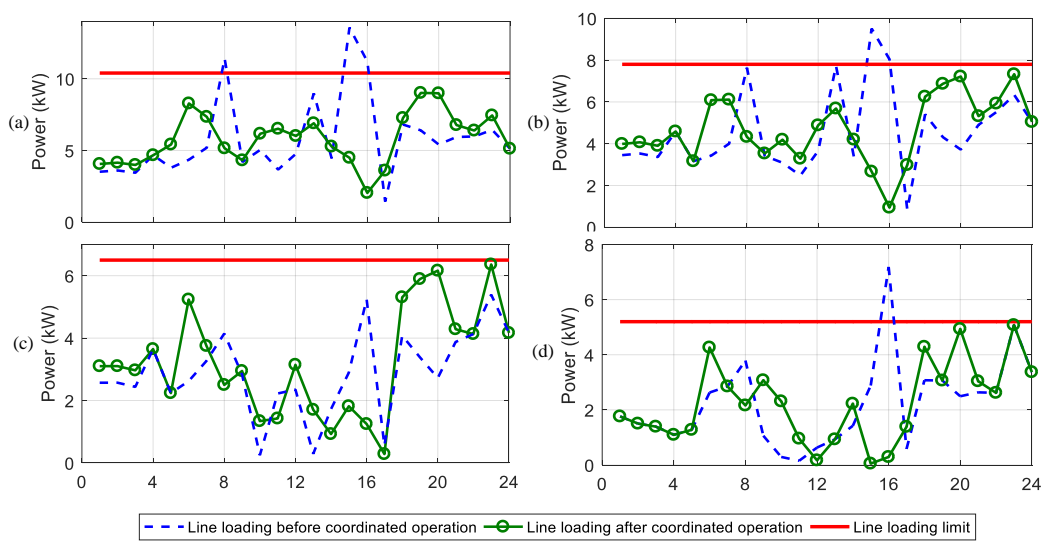

Figure 8: loading before and after implementation the proposed model for the line connecting, (a) network to the upstream grid, (b) Home #1 to #2, (c) Home #2 to #3, (d) Home #3 to #4

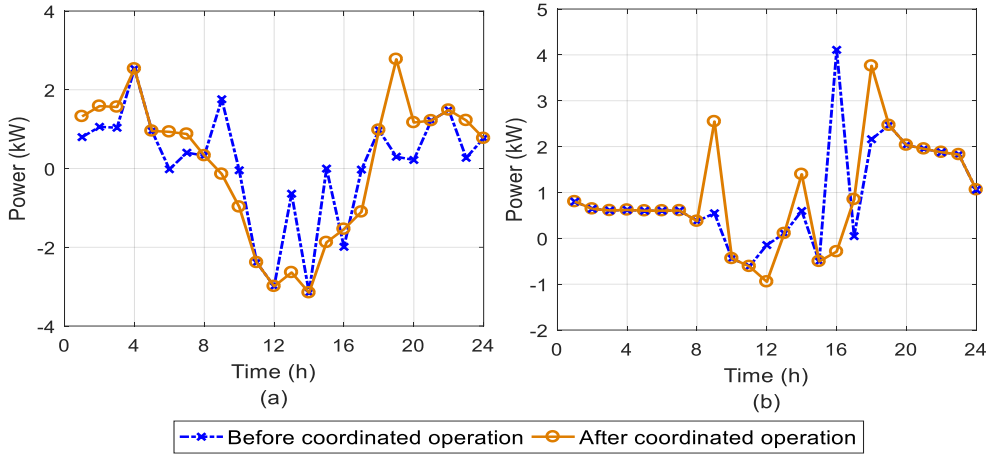

Figure 9: Consumption profile before/after implementation the proposed model for (a) Home #3, (b) Home #4

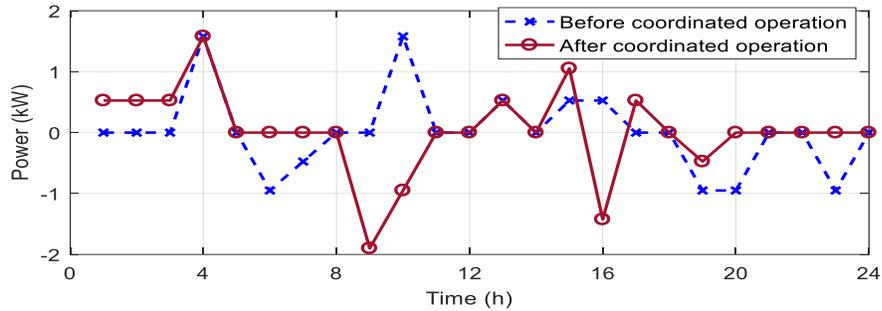

Figure 10: BESS profile before/after implementation of the proposed consumption model for Home #3

The number of iterations to reach convergence is expected to increase slightly by adding more customers as more communication might be required to reach an agreement. However, the convergence is always guaranteed because the adjustment scheme in (32) makes the performance of the algorithm less dependent on the number of optimization variables. On the other side, the scalability and computational time of the algorithm may not be a concern owning to the parallel action of customers in solving their optimization problem.

### 5.4    Impact of parameter $w$

The parameter $w$ in (23) controls the amount of reward or penalty given to the customers. The reward coefficient ($\sigma_t$)is a multiplication of $1/w$ and an exponential function in which the parameter $1/w$ is included in the power value as well. Thus, it is expected that by increasing the value of $w$, the reward value declines. By similar reasoning, the penalty coefficient likely decreases with an increase in $w$. However, as the value of $w$ increases, the rate of decreases for the reward value is much bigger than that of the penalty value. In other words, at a certain value of $w$, the penalty value tends to become larger than the reward value, as illustrated in Figure 13. Altogether, by increasing $w$, the electricity payment of customers expected to rise and accordingly the profit of the DNO is likely to grow.

The impact of different values of $w$ on the payment of customers and the profit of DNO is investigated and summarised in Figure 14. By increasing $w$, while the payment of customers increases but they still pay less bill compared to the uncoordinated HEMS operation.  As the value of $w$ reaches 0.7, the payment of customers goes above the payment they made before coordination. On the other hand, the profit of the DNO expands as $w$ increases since more payments are made by network customers. For small values of $w$, the DNO profit is reduced compared to the uncoordinated HEMS operation. However, this profit is just associated with purchasing electricity from the upstream network and the customer's payment. In other words, implementing the proposed DR program provides other savings for the DNO in terms of the damage costs for the network components which might be higher than the total rewards given to the customers to participate in the DR program. Therefore, the DNO still makes a profit for any value of $w$. Generally, the DNO can adjust the value of $w$ given his/her own financial and technical preferences and considering how much reward he/she is willing to offer to the customer to make them participate in the coordinated load scheduling program.

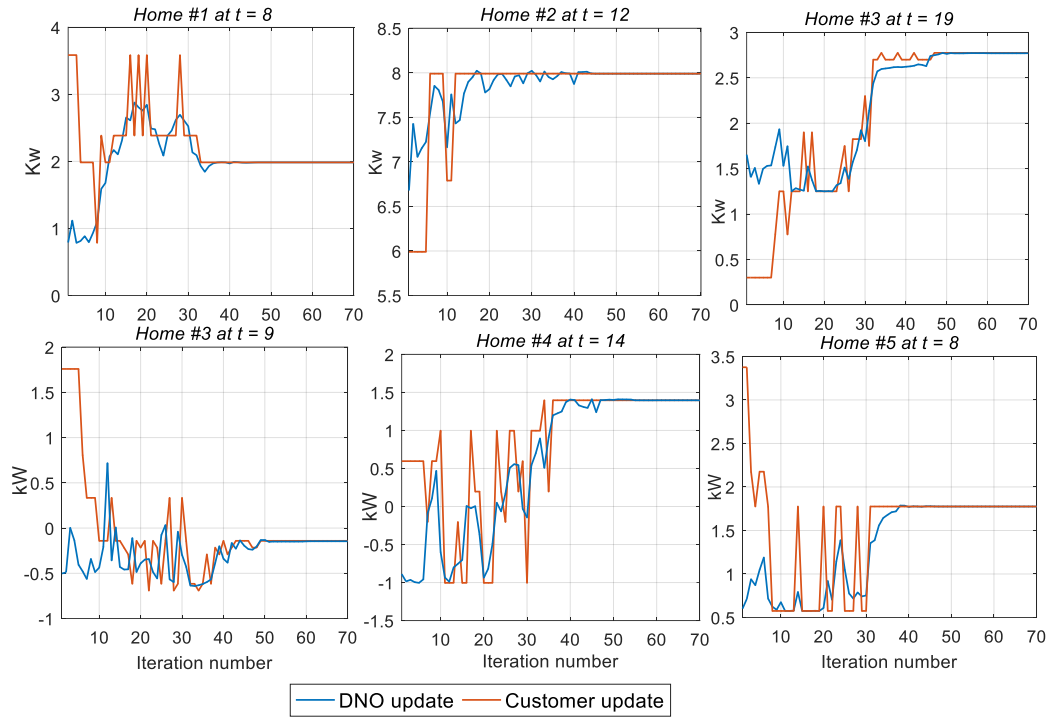

Figure 11: Load adjustment process of customers at different hours during the coordinated load scheduling

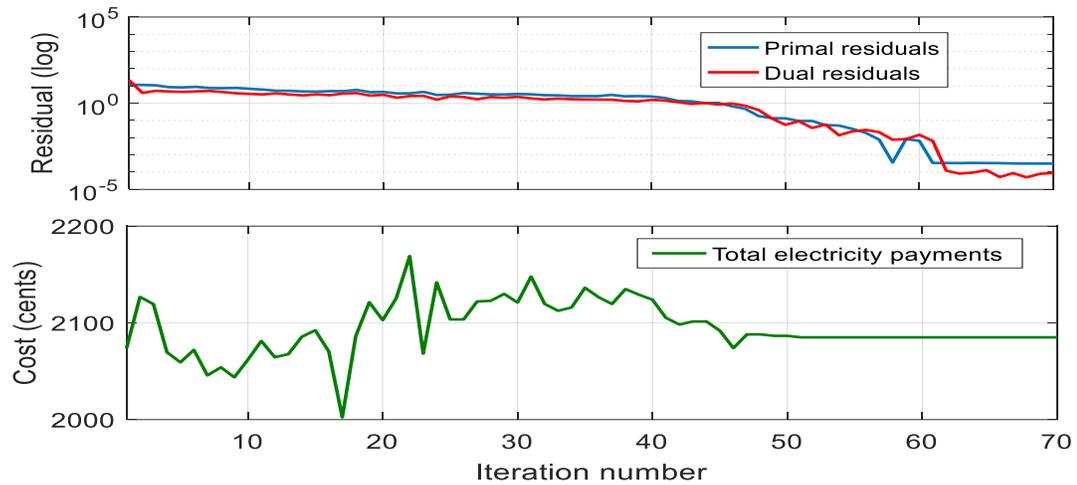

Figure 12: Convergence characteristic of the ADMM algorithm for the distributed load scheduling problem

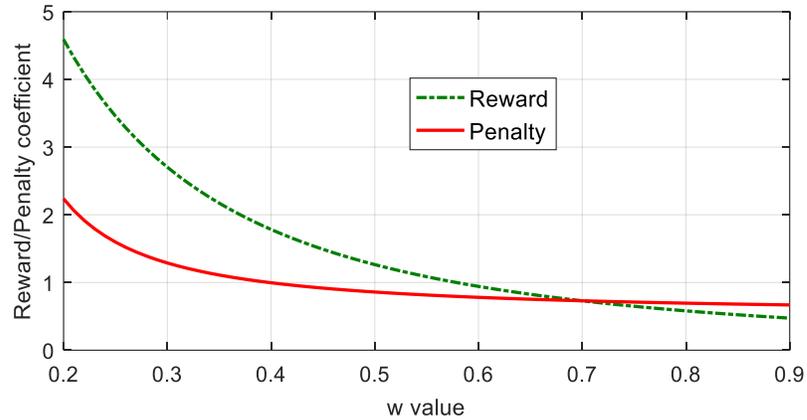

Figure 13: Reward/Penalty coefficient value as a function of the $w$ parameter

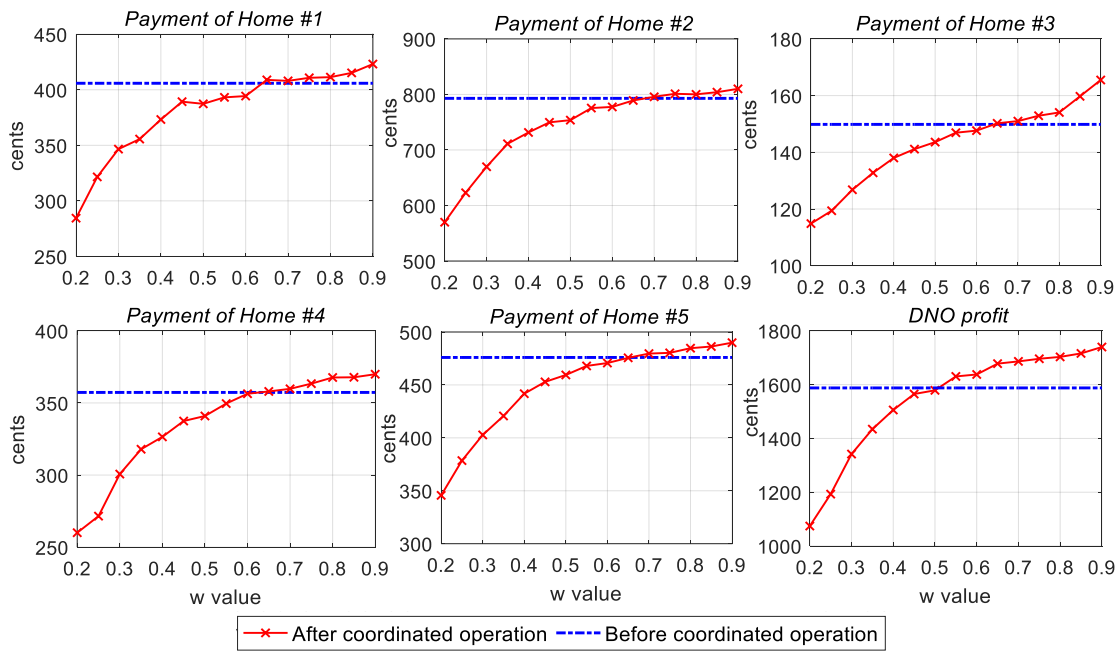

Figure 14: Effect of changing the parameter $w$ on the electricity bill of customers and the DNO profit

### 5.5    Individualized incentives

The price gap value in (21) and the reward/penalty coefficient in (23) are calculated based on the aggregated load of the network and then offered to all customers. Therefore, all customers receive the same reward or penalty regardless of their contribution to reducing the mismatch between the available supply and aggregated load. A customer may also receive a reward/penalty without any load shifting due to the change in the aggregated load of the network. In other words, some customers might be benefited from the contribution of others. For example, in Figure 9, from 1 a.m. to 8 a.m., customer #4 does not change its consumption. However, because of the load shifting of customer #3, the network load profile has improved in this period and consequently, all customers including customer #4 are benefited from the price reductions. In such a situation, some customers do not follow the DNO's suggested consumption but still get benefited.

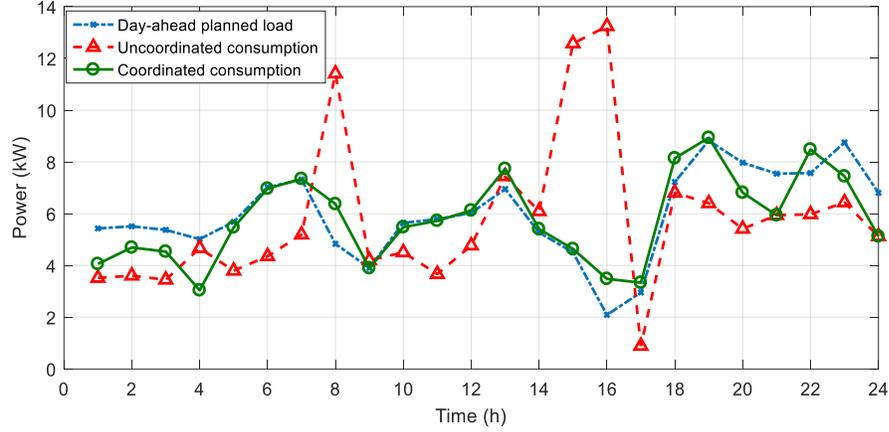

Figure 15: Network load before and after HEMS coordination based on the individualized incentives

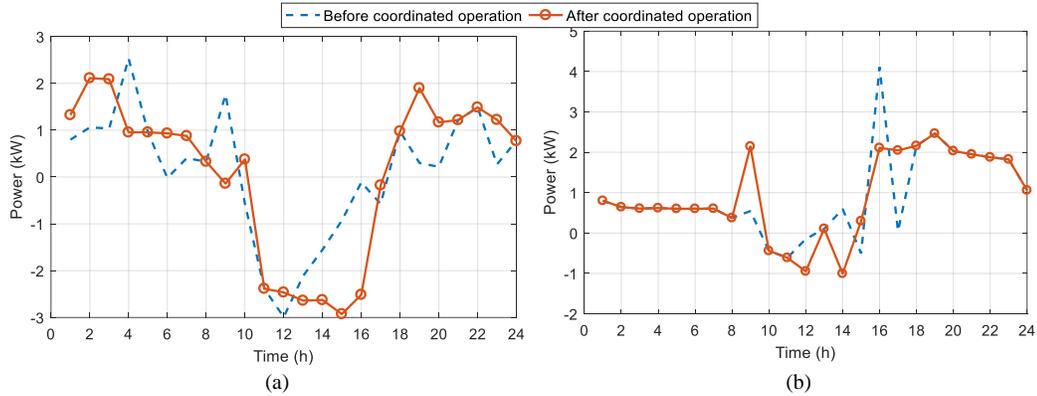

Figure 16: Consumption profiles based on the individualized incentives for (a) Home #3, (b) Home #4

One way to resolve this is by providing "*individualized*" price incentives to every customer rather than offering the same "*global*" incentives to all customers. This means that each customer receives separate reward/penalty values depending on her response to the suggested consumptions by the DNO. In this regard, the reward/penalty coefficient for each customer $h$ at time $t$ is calculated using (34).

$$\sigma_{t,h} = \begin{cases} w^{-1}\left(e^{\left(\frac{1}{w}\right)\left(\frac{\hat{L}_{h,t} - L_{h,t}^{unc}}{L_{h,t}^{unc}}\right)} - 1\right), & \Delta L_t > 0 \\ w\left(e^{\left(\frac{1}{w}\right)\left(\frac{\hat{L}_{h,t} - L_{h,t}^{unc}}{L_{h,t}^{unc}}\right)} - 1\right), & \Delta L_t < 0 \end{cases} \tag{34}$$

And in the same way, the FiT price is obtained by:

$$\sigma_{t,h} = \begin{cases} 0, & \Delta L_t > 0 \\ w^{-1}\left(e^{\left(\frac{1}{w}\right)\left(\frac{\hat{L}_{h,t} - L_{h,t}^{unc}}{L_{h,t}^{unc}}\right)} - 1\right), & \Delta L_t < 0 \end{cases} \tag{35}$$

Like the global incentives, the value of $\sigma_{t,h}$ is non-zero when $\Delta L_t \neq 0$. However, the $\sigma_{t,h}$ value depends on the contribution of each customer in minimizing supply-demand imbalance and is not a function of the aggregated load of the network anymore. Then, the price gap value for each customer $h$ at time $t$, i.e., $\in_{t,h}$, is calculated by (21). Therefore, real-time individualized prices offered to each customer obtained by:

$$\lambda_{t,h}^{RT} = \lambda_t + \in_{t,h} \qquad (36)$$

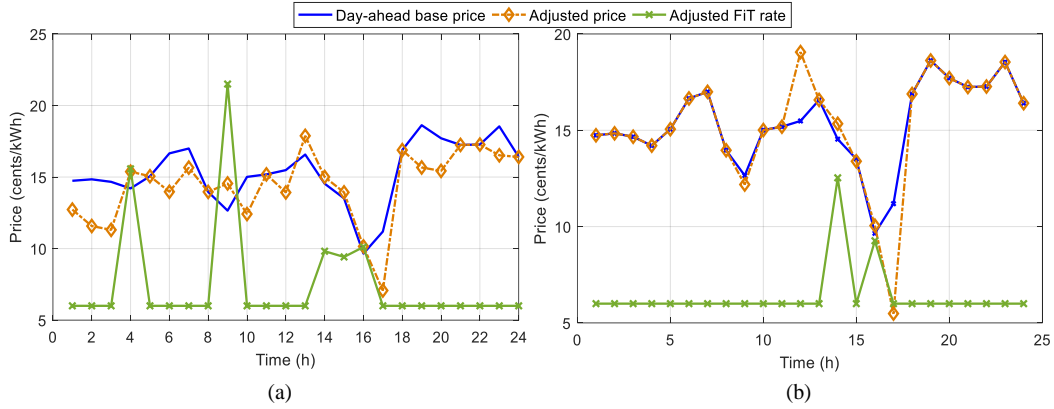

Figure 17: Adjusted price with individualized reward/penalty values added, (a) Home #3, (b) Home #4

Table 3: Payment comparison of customers in two different incentive design schemes for w=0.35

| | Electricity payment (cents) | | |
|---|---|---|---|
| | uncoordinated operation | coordinated operation with global incentives | coordinated operation with individualized incentives |
| Home #1 | 405.872 | 355.630 | 397.810 |
| Home #2 | 792.959 | 711.133 | 781.535 |
| Home #3 | 149.847 | 132.770 | 145.943 |
| Home #4 | 357.306 | 317.905 | 349.124 |
| Home #5 | 475.961 | 420.677 | 467.402 |

Then, by substituting (36) in (8) and computing the payment of each customer and then employing the auxiliary variables, the same equation as in (29) is found for the DNO update phase. The only difference is that $\in_t$ is replaced by $\in_{t,h}$.

The comparison for the aggregated load of the network after implementing the proposed model employing individualized incentives is done in Figure 15. All rebound peaks at 8 a.m. and 3-4 p.m. have been successfully removed.

The detailed load schedule for home #3 and #4 are depicted in Figure 16. The adjusted price values with individualized reward/penalty values added are displayed in Figure 17. The price curves in the real-time operation of the network are different for each home and depend on the individual contribution of customers. The electricity price and the FiT price is only adjusted at any time the customer changes its original consumption. For example, from 1 a.m. to 8 a.m., customer #4 has not changed her consumption, and accordingly, no incentive is added to the price. However, customer #3 adjusted her consumption dynamically in this period to help the DNO reduce the supply-demand mismatch in the network. Therefore, she receives lots of rewards in the form of price reduction.

The electricity payments of customers before and after coordination are compared in Table 3. The results are listed for *global* and *individualized* incentives. For both cases, electricity payments of customers are reduced compared with the uncoordinated operation case, however, more reduction happens for the global incentives. This is because, in the case of global incentives, some customers at some hours of the day are benefited from the contribution of others while in the individualized incentive case rewards are only given based on the specific contribution of customers.

## 6    Conclusions and Future Work

This paper aims to develop a coordinated load scheduling framework to prevent undesirable impacts of "rebounding peaks" due to the individual operation of multiple HEMS on the network. The objective is to design a system-wide framework for coordinating the operation of multiple HEMS to shape the aggregated load profile of the neighborhood. At the same time, the proposed approach ensures that the customers are not paying more than they would have paid under individually optimized HEMS. Additionally, it aimed to maintain the comfort degree of customers by meeting the specified scheduling deadline constraints of users. In this regard, appropriate reward and penalty values are designed and offered to the customers to incentivize them to participate in such a load scheduling program and adjust their load profile. Reward/penalty values are included in the electricity bill of the customers in the form of adjustments to the electricity prices. Two mechanisms for designing appropriate reward/penalty values are formulated and *global* and *individualized* pricing structures are proposed. Both pricing structures prove to work efficiently in keeping the operational constraints of the network within desirable limits and reducing the electricity bills of customers. Then, it is the choice of the DNO to select the proper pricing structure to adequately incentivize customers to join the load coordination framework.

The ADMM algorithm is selected to solve the formulated problem in a distributed manner. The distributed optimization model requires a minimum amount of information exchange between customers and the DNO ensuring customer privacy. The novelty of the proposed framework lies in its ability to coordinate multiple HEMS in the residential neighborhood by considering the electricity network constraints. The presented simulation results demonstrate the efficiency of the proposed method and the effectiveness of the distributed algorithm in maintaining the supply-demand balance of the network and reducing the electricity bills. Furthermore, the DNO's make more profit due to the elimination of network peak and associate damage costs.

Although the ADMM algorithm is a good choice and well-suited for the distributed energy management implementation in theory, it is typically exposed to errors due to communication noise, delay, or loss of signals and lack of enough data. Therefore, the impact of these errors on the performance of the proposed ADMM-based coordination of HEMS will be investigated in the future requiring separate mathematical modeling and formulation.

## Appendix A.    HEMS Model

Each home $h \in H$ is assumed to have $|A_h|$ number of deferrable appliances, where $A_h$ is the set of all deferrable appliances such as washing machine and dishwasher. Also, some customers might have installed solar PV or BESS. Each deferrable appliance $a \in A_h$ might have a different power consumption $P_h^a$, a different job length $l_h^a$ to finish the task and a different operation period $[\alpha_h^a, \beta_h^a]$. In the proposed HEMS, $\boldsymbol{s}_{h,t} = \left( s_{h,t}^a, s_{h,t}^b \right)$ and $\boldsymbol{u}_{h,t} = \left( u_{h,t}^a, u_{h,t}^b \right)$ represent the total state and decision variables of the problem for each customer $h$ at time $t$, respectively.

The state of each appliance $a$ for home $h$, i.e., $s_{h,t}^a$, is defined by the number of time slots that the appliance has been operated, i.e., $s_{h,t}^a \triangleq l_t^a$, where $l_t^a \in \{0,1,\ldots,I_h^a\}$ is the operational status of the appliance. In this model, $s_{h,t}^a = 0$ means that the appliance is not activated yet by the HEMS while $s_{h,t}^a = I_h^a$ shows that the operation has been completed. The decision variable $u_{h,t}^a \in \{0,1\}$ controls the operation of the appliance. When $t \notin [\alpha_h^a, \beta_h^a]$, the appliance cannot be operated and consequently $u_{h,t}^a = 0$ and $s_{h,t}^a = 0$. During the time interval $t \in [\alpha_h^a, \beta_h^a - I_h^a + 1)$, the state of the appliance changes by $s_{h,t+1}^a = s_{h,t}^a + u_{h,t}^a$. At $t = \beta_h^a - I_h^a + 1$, if $s_{h,t}^a = 0$, it should be activated immediately. Once activated (i.e., $u_{h,t}^a = 1$), the appliance has to work for $I_h^a$ time slots without any interruption until the job is done in which case the state moves to $s_{h,t}^a = I_h^a$ and the appliance has to remain idle for the rest of the day.

For the BESS, the state is represented by the state of charge (SOC), i.e., $s_{h,t}^b \triangleq SOC_t$. The SOC is constrained by $SOC_t \in [SOC_{min}, SOC_{max}]$. The decision variable for the operation of BESS is denoted by $u_{h,t}^b$ which handles charging (when $u_{h,t}^b \geq 0$) and discharging (when $u_{h,t}^b \leq 0$) process and is constrained by charging/discharging limits, i.e., $u_{h,t}^b \in [\Delta SOC_{min}, \Delta SOC_{max}]$. The evolution of the BESS state is governed by $s_{h,t+1}^b = s_{h,t}^b + u_{h,t}^b$. The BESS power of home $h$ at each time $t$, $P_{h,t}^b$, is given by:

$$P_{h,t}^b = \frac{E_h^b}{\eta_c^b \Delta t}(u_{h,t}^b)^+ + \frac{E_h^b \eta_d^b}{\Delta t}(u_{h,t}^b)^- \tag{A.1}$$

where $(.)^+ = max\{0, .\}$ and $(.)^- = min\{0, .\}$. $E_h^b$ is the capacity of the BESS at home $h$ and $\Delta t$ is the duration of each time slot. $\eta_c^b$ and $\eta_d^b$ are the charging and discharging efficiencies, respectively. Positive value for $P_{h,t}^b$ means charging and negative $P_{h,t}^b$ indicates discharging of the BESS.

The total electricity consumption of home $h$ at time $t$ is calculated by:

$$L_{h,t} = \sum_{a \in A_h} P_h^a u_{h,t}^a + P_{h,t}^{non} + P_{h,t}^b + P_{h,t}^{pv} \tag{A.2}$$

where $P_{h,t}^{pv}$ and $P_{h,t}^{non}$ are the PV generation and non-controllable demand of home, respectively. The non-controllable demand includes the basic needs of the household such as lighting, TV, and refrigerator. The HEMS has no control over this load type, and it has to be met at any time without any interruption.

The goal of HEMS is to decide when to turn on and turn off home appliances and find the best operational mode of battery to minimize the customer's electricity payment. To this end, each HEMS solves an optimization problem with the total objective function of:

$$J_h(L_{h,t}) = \min_{\substack{u_{h,t}^a \in \Phi_{h,t}^a \\ u_{h,t}^b \in \Gamma_{h,t}^b}} \sum_{t=1}^{T} C_{h,t}(s_{h,t}, u_{h,t}) \tag{A.3}$$

where $C_{h,t}(s_{h,t}, u_{h,t})$ models the one-stage cost or the cost of being at state $s_{h,t}$ and making decision $u_{h,t}$ in the proposed DDP model. respectively. $\Phi_{h,t}^a$ is the set of all feasible operating decisions for each controllable appliance of home $h$ at stage $t$ based on the related constraints. The set of feasible decisions for the BESS defined by $\Gamma_{h,t}^b$ as follows:

$$\Gamma_{h,t}^b = \{u_{h,t}^b | u_{h,t}^b \in [\Delta SOC_{min}, \Delta SOC_{max}], \ SOC_{min} \leq SOC_t \leq SOC_{max}\}$$

Problem (A.3) is a multi-stage sequential-decision problem and the optimal solution can be found using the algorithm developed in [31].

## Appendix B.    Preliminaries of ADMM

The alternating direction method of multipliers which is abbreviated as ADMM is a technique developed for solving large scale and complex optimization problems. It divides the original problem into several small-scale sub-problems which are easy to solve individually. Then, the solution of the original problem is obtained by solving these sub-problems in a parallel manner. In its general form, ADMM solves a problem in such form [28]:

$$\min_{x,z} f(x) + g(z)$$
$$subject\ to:\ Ax + Bz = c \tag{B.1}$$

where $x \in \mathbb{R}^n$ and $z \in \mathbb{R}^m$ are optimization variables, $c \in \mathbb{R}^p$, $A \in \mathbb{R}^{p \times n}$ and $B \in \mathbb{R}^{p \times m}$ are matrices containing known parameters. $f(.)$ and $g(.)$ model two convex objective functions. Please note that two objective functions have separate optimization variables and may have their own constraints as well, but an equality constraint couples them together.

The scaled augmented Lagrangian of (B.1) is described by:

$$\mathcal{L}_\rho(x,y,u) = f(x) + g(z) + \frac{\rho}{2}\|Ax + Bz - c + u\|_2^2 \tag{B.2}$$

where $\rho > 0$ is the penalty parameter and $u$ is the vector of scaled dual variables. $\|.\|_2$ represents the $l_2$-norm of a vector.

Under the ADMM framework, $x$ nad $z$, as the optimization variables, are obtained and updated in an alternating and sequential fashion iteratively using the scaled dual variables $u$ as the link between the iterations and optimization variables. In this method, every iteration consists of a $x$-update step in which $x$ is solved considering $z$ as a fixed value. Then, in the $z$-update step, the problem is solved for $z$ with $x$ as a fixed value as found previously in the $x$-update step. Finally, scaled dual variables $u$ are updated for the next iteration. Mathematical description of each iteration $k + 1$ of ADMM is as follows [28]:

$$x^{k+1} = \underset{x}{\operatorname{argmin}} \ f(x) + \frac{\rho}{2}\|Ax + Bz^k - c + u^k\|_2^2 \tag{B.3}$$

$$z^{k+1} = \underset{z}{\operatorname{argmin}} \ g(x) + \frac{\rho}{2}\|Ax^{k+1} + Bz - c + u^k\|_2^2 \tag{B.4}$$

$$u^{k+1} = u^k + Ax^{k+1} + Bz^{k+1} - c \tag{B.5}$$

The convergence behavior of the ADMM has already been studied and the proof of convergence is obtained. It has been shown that this method converges to the optimal solution under some assumptions [28]. So, explaining the proof of convergence is not within the scope of this research.

For the ADMM technique, the termination criteria are usually defined by:

$$\|r^k\|_2 = \|Ax + Bz - c + u^k\|_2 \leq \epsilon^{prim} \tag{B.6}$$

$$\|s^k\|_2 = \|\rho A^T B(z^k - z^{k-1})\|_2 \leq \epsilon^{dual} \tag{B.7}$$

where $r^k$ is the primal residual, $s^k$ the dual residual of the algorithm and $\epsilon^{prim} > 0$ and $\epsilon^{dual} > 0$ are primal and dual tolerances parameters, respectively. The algorithm stops when tolerant parameters are adequately small. A reasonable value might be $\epsilon^{prim} = \epsilon^{dual} = 10^{-3}$ or $10^{-4}$.

## Appendix C.   Test feeder data and simulation parameters values

| | Parameter | Unit | Value |
|---|---|---|---|
| Feeder data* | $R$ | % on 1 MVA base | 153 |
| | $X$ | % on 1 MVA base | 62.5 |

*All lines have the same R nd X parameters

| | Parameter | Unit | Value |
|---|---|---|---|
| | $T$ | hr | 24 |
| | $\Delta t$ | hr | 1 |
| | $a$ | $cents/kWh^2$ | 0.2 |
| | $b$ | $cents/kWh$ | 2 |
| Simulation parameters | $c$ | - | 0 |
| | $P_t^{DG}$ | kW | 0 |
| | $\varphi$ | - | 4.8 |
| | $FiT_t$ | cents/kWh | 6 |
| | $\mu_t^b$ | cents/kWh | 2 |
| | $\mu_t^s$ | cents/kWh | 2 |